\documentclass[11pt]{article}

\usepackage{times}

\usepackage{amssymb}
\usepackage{amsmath}
\usepackage{amsthm}
\usepackage{amscd}

\usepackage[xdvi]{graphicx}
\usepackage{subfigure}

\newtheorem{theorem}{Theorem}
\newtheorem*{theorem*}{Theorem}
\newtheorem{lemma}{Lemma}
\newtheorem{corollary}{Corollary}
\newtheorem{proposition}{Proposition}

\newcommand{\Tn}{\ensuremath{\mathcal{T}_{n}}}

\newcommand{\PT}[1]{\ensuremath{p_{#1}}}
\newcommand{\asymPT}[1]{\ensuremath{p'_{#1}}}
\newcommand{\rootPT}[1]{\ensuremath{r_{#1}}}
\newcommand{\asymrootPT}[1]{\ensuremath{r'_{#1}}}
\newcommand{\NC}{\ensuremath{{\rm NC}(n)}}
\newcommand{\enorb}[2]{\ensuremath{o(#1, #2)}}

\newcommand{\dist}[2]{\ensuremath{\delta(#1,#2)}}

\newcommand{\tp}[1]{\ensuremath{U_{#1}}}
\newcommand{\bt}[1]{\ensuremath{L_{#1}}}

\newcommand{\aarc}{\ensuremath{\rightharpoonup}}
\newcommand{\barc}{\ensuremath{\rightharpoondown}}

\newcommand{\mP}{\ensuremath{\mathcal{P}}}

\newcommand{\N}{\ensuremath{\mathbb{N}}}

\begin{document}

\title{\vspace*{-2ex}Counting orbits under Kreweras complementation}

\author{Christine Heitsch}

\date{Georgia Institute of Technology \\
School of Mathematics \\
Atlanta, GA 30332 - 0160 \\
\texttt{heitsch@math.gatech.edu}}

\maketitle

\begin{abstract}

The Kreweras complementation map is an anti-isomorphism
on the lattice of noncrossing partitions.
We consider an analogous operation for plane trees motivated by the 
molecular biology problem of RNA folding. 
In this context, we explicitly count 
the orbits of Kreweras' map according to their length as the number of
appropriate symmetry classes of trees in the plane.
These enumeration results are consolidated into a single implicit 
formula under the cyclic sieving phenomenon. 
\\


\end{abstract}

\section{Noncrossing partitions}\label{nonpart}

For $i, j \in \N$, let $[i,j] = \{k \in \N \mid i \leq k \leq j\}$.
Let $\mP$ be a (set) partition of $[1,n]$.
Then $\mP$ is \emph{noncrossing} if there do not exist 
$a, c \in P$ and $b,d \in Q$ for distinct parts $P, Q \in \mP$
such that $a < b < c < d$.
Noncrossing partitions are the subject of much 
study~\cite{mccammond-06, simion-00},
both for their intrinsic interest as well as connections to other 
fields.

Let \NC\ denote the set of noncrossing partitions ordered
under refinement.
This is a lattice, with the trivial partition 
$\{[1,n]\}$ as the top element and the singleton 
one $\{\{i\} \mid 1 \leq i \leq n\}$ as the bottom.
Moreover, it is complemented, i.e.\ for each $\mP$ 
there exists some $\mP'$ such that their least upper bound
is the top element and greatest lower bound the bottom one.
Such a $\mP'$ is explicitly given by the Kreweras complementation
map~\cite{kreweras-72}, denoted here $\kappa$.

To compute $\mP' = \kappa(\mP)$, 
interleave a second copy of $N = [1,n]$ 
with the first as $1 < 1' < 2 \ldots n < n'$.
Fix $\mP$ on $N$.
Then $\kappa(\mP)$ is the coarsest partition from $\NC$ on 
$N' = \{1', 2', \ldots, n'\}$ such that $\mP \cup \mP'$ is 
a noncrossing partition on the totally ordered set $N \cup N'$.
Geometrically, the $2n$ integers can be visualized as points on a 
circle, ordered clockwise. 
Each $P \in \mP$ is drawn as the convex hull containing 
all $i \in P$.
Then $\mP$ is noncrossing exactly when no two hulls intersect.
The image of $\mP$ under $\kappa$ is obtained by ``filling in''
the largest convex hulls on $N'$ such that there are no intersections.

The map $\kappa$ is a lattice anti-isomorphism, but  
not an involution for $n > 2$.
In fact, there is a unique complement for \mP\ exactly when 
it is either the top or bottom element.
Since $\kappa^{2n}(\mP) = \mP$,
iterating yields an orbit $\mathcal{O}$ of length $2 \leq l \leq 2n$.
The size and structure of these orbits is of
interest~\cite{simion-ullman-91}, 
in part for their connection to free probability~\cite{nica-speicher-06},
and will be characterized through an analogous operation on plane trees.

Although $\kappa$ is not an involution, 
$\kappa^{2i}(\mP)$ is a ``rotation'' of $\mP$
in the sense that $j \in P \in \mP$
exactly when $j' \equiv j - i \pmod n \in P' \in \kappa^{2i}(\mP)$.
In contrast, the block structure of $\mP$ and $\kappa(\mP)$
is typically quite different.
However, we show that $\kappa(\mP)$ is also a ``rotation'' 
of \mP\ since they both correspond to the same plane tree.
Hence, 
the problem of counting orbits under Kreweras complementation 
reduces to known solutions for enumerating different symmetry 
classes of trees in the plane.

\section{Complementing plane trees}\label{trees}

A \emph{plane tree} is a rooted tree whose subtrees are linearly 
ordered~\cite{stanley-99}.
They are also known as ordered or linear trees and, 
like noncrossing partitions, are enumerated by the Catalan numbers,
$C_{n} = \frac{1}{n+1}\binom{2n}{n}$.  
Let \Tn\ denote the set of plane trees with $n$ edges, and $T \in \Tn$.
Motivated by the molecular biology problem of RNA folding~\cite{barrier}, 
label the boundary of $T$ 
with $[1,2n]$ in increasing order counter-clockwise from the root.
Let $(i,j)$ denote the edge in $T$ with indices $i < j$ on the 
left and right sides.
There is a straight-forward bijection between plane trees expressed
as $T = \{(i,j) \mid 1 \leq i < j \leq 2n\}$ and noncrossing
perfect matchings on $2n$ endpoints.

In a rooted tree $T$, the degree of a vertex 
is the number of its children.
A plane tree can be defined recursively as
$T = (r, T_{1}, T_{2}, \ldots, T_{k})$
where root vertex $r$ is connected to the roots of $k$
plane trees in linear order.
Given an unrooted, unordered tree embedded in the plane,
the choice of root is determined both by $r$ and by
selecting a first child $T_{1}$ among the cyclically ordered 
subtrees $T_{1}, T_{2}, \ldots, T_{k}$.
 
Let $\phi(T) = \{(k-1,2n) \mid (1,k) \in T\} 
\cup \{(i-1,j-1) \mid (i,j) \in T,\ 1 < i < j \leq 2n\}$.
Then $\phi$ is a bijection on \Tn\ whose action
corresponds to rerooting $T$
by shifting the indices one edge counterclockwise.
We claim this corresponds to Kreweras complementation.

There are various bijections between \NC\ and \Tn.
Classically~\cite{dershowitz-zaks-86, prodinger-83},
a partition with $k$ parts corresponds to a tree with $k$ leaves.
More recently~\cite{bona-etal-00, callan-smiley-05, liu-wang-li-14},
ones for bicolored plane trees have been given.
Here we describe a map which is a special case of the general 
``tree partition'' definition~\cite{barrier}.
This bijection is dual to~\cite{bona-etal-00, callan-smiley-05}, 
and this duality illuminates the relationship between $\mP$ and $\kappa(\mP)$.

First, observe that exactly one index for each edge is odd.
Call $(i,j)$ \emph{odd} if $i$ is and \emph{even} otherwise.
Edges incident on the root are all odd, and edge parity
alternates along a path from the root to a leaf.
Call the root an even vertex.
Otherwise, a vertex inherits its parity from the 
incoming, i.e.\ parent, edge.

The bijection $\rho: \Tn \rightarrow \NC$ is obtained from a 
tree partition for $T$ and 
$\bt{n} = \{(1,2n)\} \cup \{(2i, 2i+1) \mid 1 \leq i \leq n - 1\}$.
In this case, 
the parts of the tree partition correspond to subtrees of $T$ consisting
of all the odd children of a vertex, and their even parent (if there
is one).
A noncrossing partition, c.f.\ Theorem~8 of~\cite{barrier},
is obtained by projecting the odd indices 
for each subtree down as $2i - 1 \rightarrow i$.

Dually, $\tp{n} = \{(2i-1, 2i) \mid 1 \leq i \leq n\}$ could be used,
which recapitulates~\cite{bona-etal-00, callan-smiley-05}.
Now, however, the subtrees of $T$ corresponding to the tree partition 
with $\tp{n}$ consist 
of all the even children of a vertex, and their odd parent.
Let $\bar{\rho}$ denote the bijection with \NC\ obtained by 
projecting the odd indices for each subtree down to a subset of $[1,n]$.

Observe that $\phi(\tp{n}) = \bt{n}$, and vice versa.
Both are ``star'' trees with $n$ edges that differ only in the choice
of root.
They are also the top and bottom element in the \Tn\ poset~\cite{barrier}
under the bijection $\rho$.
Moreover, their corresponding noncrossing perfect matchings form 
a meander~\cite{meander}, 
i.e.\ a single closed loop when drawn on the same set of endpoints,
one set above and the other below the line.

More generally, let $T' \in \Tn$.  
Then $T$ and $T'$ are complements in the lattice sense exactly 
when~\cite{barrier, meander}
their corresponding noncrossing perfect matchings form a meander.
Let $i \aarc j$ denote the edge in $T$ with odd index $i$, and
$j \barc i$ the one in $T'$.
Then $T$ and $T'$ are complements if there is a single closed loop
consisting of the $2n$ distinct indices
starting from $(1,j) \in T$ and finishing with $(1, j') \in T'$;
\[ 1 \aarc j \barc i \ldots i' \aarc j' \barc 1\text{.} \]

\begin{lemma}
$T$ and $\phi(T)$ are complements.
\end{lemma}

\begin{proof}
Let $(1,2k) \in T$.  Then $(2k-1, 2n) \in \phi(T)$.
Let $T'$ be the subtree of T with indices $[2,2k-1]$, and $T''$ the one
for $[2k+1,2n]$.
Inductively, they form meanders with their images under $\phi$.
Keeping the indices from $T$, we have single closed loops $L'$ and $L''$
respectively;
\[ 
2 \aarc j' \ldots j' - 1 \barc 2k -1 \aarc i' \ldots \aarc 2\ \text{and}\ 
2k+1 \aarc j'' \ldots i'' \aarc j'' - 1 \barc 2n \ldots \aarc 2k+1\text{.} \]
However, in $\phi(T)$, the edge $(2,j')$ is mapped to $(1, j'-1)$ 
rather than $(j'-1, 2k-1)$ as in $L'$.
Likewise, for $(2k+1, j'') \in T$,
there would be $j'' - 1 \barc 2k$ for $\phi(T)$ rather than 
$j'' - 1 \barc 2n$ as in $L''$.
But then, moving backwards through $L''$ and forward through $L'$, 
\[ 1 \aarc 2k \barc j'' - 1 \aarc i'' \ldots j'' \aarc 2k+1 \ldots 2n 
\barc 2k - 1 \aarc i' \ldots 2 \aarc j' \ldots j'-1 \barc 1\]
is a single closed loop.
\end{proof}

Let $(i,j), (i', j') \in T$.
Note that $(i',j')$ is the first child of $(i,j)$ if $i' = i+1$,
respectively last if $j' + 1 = j$.
Likewise, $(i', j')$ is the next sibling of $(i,j)$ if $i' = j+1$,
or previous if $j' + 1 = i$.

\begin{theorem}\label{commute}
$\rho(\phi(T)) = \kappa(\rho(T))$.
\end{theorem}

\begin{proof}
Consider the noncrossing perfect matching for $T$, denoted $C$, 
as $n$ chords on $2n$ labeled endpoints clockwise around a circle.
For $1 \leq i \leq n$,
let $a_i$ denote the arc with endpoints $2i - 2 \pmod{2n}$ and $2i - 1$,
and $b_i$ the one with $2i-1$ and $2i$.
Note the edges in \bt{n}\ correspond to $a_i$, and those in 
\tp{n}\ to $b_i$.

Let $\mP = \rho(T)$.
Let $S$ be a subtree of $T$ consisting of all the odd children of a
vertex, and their even parent if there is one.
Then $S$ corresponds to a single closed loop on $C$ consisting of
the edges from $S$, alternating with $a_i$,
in linear order of the children, and finishing with the parent.
If there is no parent, the first child is $(1,2k)$ and 
$(2k'-1, 2n)$ the last.
Then the image of $S$ under $\rho$ are the indices of the $a_i$.

Dually, let $\mP' = \bar{\rho}(T)$.
Let $S'$ be a subtree of $T$ consisting of all the even children of a
vertex, and their odd parent.
Then, as before, the edges from $S'$, alternating with $b_i$, 
form a single closed loop, and
the image of $S'$ under $\bar{\rho}$ are the $b_i$ indices.

Observe that the loops for $\rho(T)$ and for $\bar{\rho}(T)$
trace out complementary regions of the circle's area.
Identifying the endpoints of $a_i$ yields the convex hull 
for $\rho(T)$, and likewise for $b_i$ and $\bar{\rho}(T)$.
Hence, $\bar{\rho}(T) = \kappa(\rho(T))$.
However, consider instead rotating $C$ one arc counter-clockwise.
Then $\bar{\rho}(T) = \rho(\phi(T))$.
\end{proof}

\noindent
Hence, 
the interleaved convex hulls for $\mP$ and $\kappa(\mP)$ 
correspond to a rotation, by one index, of the same plane tree.
Since the action of $\phi$ changes the root,
but leaves the rest of the structure unaltered,
the set $\phi^{i}(T)$ corresponds to one planar tree.

\section{Counting planar trees}\label{class}

Here, ``plane tree'' refers to a rooted, linearly ordered tree.
However,
it can also mean a tree which is embedded in the plane.
Alternatively, the second type is called a \emph{planar tree}.
As an aside, a rooted tree is called planted if the 
degree of the root vertex is 1.
Planted planar trees with $n$ edges, which are also planted plane trees,
are in bijection with plane trees with $n - 1$ 
edges~\cite{harary-prins-tutte-64, klarner-70}. 

Let $E$ be a planar tree with $n$ edges.
If $E$ has no nontrivial rotational symmetry,
call it \emph{asymmetric}.
Suppose $E$ is rooted at vertex $r$.
Then the subtrees of $r$ are cyclically ordered, whereas the subtrees
of all other vertices in $E$ are linearly ordered.
In this case, $E$ is asymmetric if it has no nontrivial rotational 
symmetry which preserves the root.

Let $\PT{n}$ denote the number of planar trees with $n$ edges, 
and $\asymPT{n}$ the asymmetric ones.
Let $\rootPT{n}$ denote those that are rooted, and 
$\asymrootPT{n}$ the asymmetric rooted ones.
Then there are explicit enumeration 
formula~\cite{bergeron-labelle-leroux-98,harary-prins-tutte-64,walkup-72}
in terms of 
Euler's function $\phi(m)$,
the M\"obius function $\mu(m)$,
and the characteristic function for odd integers $\chi_{\mbox{odd}}(m)$
for $m \in \N$.
In particular,
\begin{equation}
\PT{n} = \rootPT{n} -
\frac{1}{2} C_{n} +
\frac{1}{2} \chi_{\mbox{odd}}(n) C_{\frac{n-1}{2}}\text{,}
\label{enum-PT}
\end{equation}
\begin{equation}
\asymPT{n} = \asymrootPT{n} -
\frac{1}{2} C_{n} -
\frac{1}{2} \chi_{\mbox{odd}}(n) C_{\frac{n-1}{2}}\text{,}
\label{enum-asymPT}
\end{equation}
\begin{equation}
\rootPT{n} = \frac{1}{2n} \sum_{d | n} \phi(\frac{n}{d}) \binom{2d}{d}\text{,}
\label{enum-rootPT}
\end{equation}
\begin{equation}
\text{and}\ \asymrootPT{n} = \frac{1}{2n} \sum_{d | n} \mu(\frac{n}{d}) \binom{2d}{d}\text{.}
\label{enum-asymrootPT}
\end{equation}

\section{Enumerating orbits}\label{orbits}

Let $\enorb{n}{l}$ denote the number of orbits of length $2 \leq l \leq 2n$
under $\kappa$.
Let $\mathcal{O}$ be the orbit for $\mP \in \NC$ 
and $T \in \Tn$ such that $\rho(T) = \mP$.
Since $\rho(\phi(T)) = \kappa(\rho(T))$ by Theorem~\ref{commute},
and the action of $\phi$ corresponds to rerooting $T$,
we have the following.

\begin{corollary}\label{totenum}
$\sum \enorb{n}{l} = \PT{n}$.
\end{corollary}

\begin{proof}
As an ordered, rooted tree, a plane tree is obtained from a planar tree 
by choosing a root vertex and a linear ordering, i.e.\ first child, 
for the cyclically ordered subtrees of the root.
\end{proof}

\begin{lemma}\label{exlen}
There is an orbit of length $l$ under $\phi$
only if $l = 2d$ where $d \mid n$ or $l = n$ when $n$ is odd.
\end{lemma}

\begin{proof}
Since $l \mid 2n$, suppose $l \mid n$ is odd and $l \neq n$.
For $(i,j) \in T$, let $\dist{i}{j} = \min\{j-i, 2n - j + i\}$.
Without loss of generality, suppose $(k, 2n) \in T$ has maximal
$\dist{k}{2n} = 2n - k \leq k$.
Note that $\dist{k}{2n} = n$ only if $n$ is odd and $l = n$ or $2n$.
Otherwise, $\dist{k}{2n} < n < k$.
We claim $\dist{k}{2n} < l$.
If not, $k - l < k \leq 2n - l < 2n$ which contradicts 
$(k-l, 2n-l), (k, 2n) \in T$.
Let $(i,j) \in T$.
Then, by choice of $(k,2n)$,  
$2n - l < i < k$ if and only if $2n -l < j < k$.
Hence, $k - 2n - l - 1$ must be even.
Contradiction since $k$ is also odd.
\end{proof}

\begin{theorem}\label{orbenum}
For $n > 1$, 
\[ \enorb{n}{l} = \begin{cases}
\asymPT{n} & \text{when}\ l = 2n, \\ 
C_{\frac{n-1}{2}} & \text{if}\ l = n\ \text{and $n$ is odd}, \\
\asymrootPT{\frac{n}{2}} & \text{if}\ l = n\ \text{and $n$ is even}, \\
\asymrootPT{d} & \text{for}\ l = 2d\ \text{with}\ d \mid n,\ 1 \leq d < n/2\mbox{.} 
\end{cases} \]
\end{theorem}

\begin{proof}
The unrooted, unordered planar tree $E$ corresponding to $\mathcal{O}$
is equivalent to $n$ nonintersecting cords on $2n$ unlabeled
endpoints on a circle in the same way that $T = \rho^{-1}(P)$ 
is equivalent to a noncrossing perfect matching on indices $[1,2n]$.
Hence, $l = 2n$ if and only if 
$E$ has no nontrivial rotational symmetries.
But then $\enorb{n}{2n}$ is counted by Equation~\ref{enum-asymPT}.

Similarly, $l = n$ exactly when $E$ has one rotational symmetry 
of order 2.
Consider an axis of symmetry for the circle
which splits the $2n$ points in half.
Observe that chords which cross the axis come in symmetric pairs,
unless there is a single one which is fixed.
In this case, $n$ must be odd, and we consider $(1,n-1) \in T$.
Hence, there are $C_{\frac{n-1}{2}}$ possible plane trees 
on indices $[2,n-2]$ such that the resulting $T$
have distinct orbits under $\phi$.

Suppose instead there is no cord fixed by the symmetry.
Then $n$ must be even.
Consider $(i,j) \in T$. 
Then $\dist{i}{j} < n$. 
Otherwise, the corresponding cord and its symmetric pair would 
intersect.
But then there exists some $k$ such that $\phi^{k}(T)$ consists
of two subtrees on indices $[1,n]$ and $[n+1,2n]$.

Observe that two plane trees on $\frac{n}{2}$ edges which 
differ only in the linear ordering of the root's children yield
the same set of chords for $E$.
Moreover, additional symmetries on $E$, beyond the one of order 2,
must be avoided.
Hence, when $n$ is even, $\enorb{n}{n}$ is counted by 
Equation~\ref{enum-asymrootPT} for $\frac{n}{2}$.

Finally, let $1 \leq d < \frac{n}{2}$ be an integer where $d \mid n$
and consider the circle divided into $\frac{2n}{2d}$ sectors.
By the same type of argument just used, the number of orbits
corresponding to $2d$ distinct plane trees under $\phi$ will be the
number of asymmetric rooted planar trees with $d$ edges.
\end{proof}

\section{The cyclic sieving phenomenon}\label{cyclic}

The explicit enumeration depends crucially on different rotational 
symmetries in the plane.
This is unified by counting the orbits implicitly under
the cyclic sieving phenomenon (CSP)~\cite{reiner-stanton-white-04}.
Let $X$ be a finite set and $C$ a cyclic group of order $n$ acting on $X$.
Suppose there exists a $q$-enumerator for $X$, that is a polynomial 
$X(q)$ with nonnegative integer coefficients such that $X(1) = |X|$.  
The triple $(X, X(q), C)$ exhibits CSP if either of the following
two equivalent conditions holds.
\begin{enumerate}
\item For every $c \in C$, 
\( [X(q)]_{q = \omega(c)} = |\{x \in X: c(x) = x\}|\)
where $\omega: C \rightarrow \mathbb{C}^{\times}$ is an embedding of 
$C$ into the multiplicative group of nonzero complex numbers.
\item The coefficient $a_{l}$ defined uniquely by the expansion
\( X(q) \equiv \sum_{l=0}^{n-1} a_{l} q^{l} \bmod q^{n} - 1 \)
counts the number of $C$-orbits on $X$ for which the stabilizer-order
divides $l$.
\end{enumerate}
 
\noindent
Since rooted planar trees exhibit CSP, Proposition 4.1 
of~\cite{reiner-stanton-white-04} clarifies the relationship between
\rootPT{n} and \asymrootPT{n} based on the Ramanujan sum and 
M\"obius inversion. 
Likewise, planar trees also exhibit CSP. 
\begin{proposition}~\cite{bessis-reiner-11, white-05}
Let $X$ be the set of $n$ nonintersecting chords on $2n$ labeled endpoints
around a circle,
and $C = \mathbb{Z} / 2n \mathbb{Z}$ act on $X$ by rotation. 
Let $$X(q) = \frac{1}{[n+1]_q} {2n \brack n}_q$$ be the $q$-Catalan number. 
Then $(X, X(q), C)$ exhibits the cyclic sieving phenomenon.
\end{proposition}
\noindent
The $q$-Catalan numbers are a generalization 
where the usual integers, factorials, and binomial coefficients
are replaced by their $q$-analogues starting from 
$[n]_{q} := 1 + q + q^{2} + \ldots q^{n-1}$.

For example, when $n = 4$,
\( X(q) =  1 + q^{2} + q^{3} + 2 q^{4} + q^{5} + 2 q^{6} + q^{7} + 2 q^{8} + q^{9} + q^{10} + q^{12} \)
and
\( X(q) \bmod q^{8} - 1 = 3 + q^{1} + 2 q^{2} + q^{3} + 3 q^{4} + q^{5} + 2 q^{6} + q^{7} \text{.}\) 
Thus, via CSP, there are 3 planar trees with 4 edges, 
of which one is asymmetric, another has a symmetry of order 2, and the 
third of order 4.
Alternatively,  
${\rm NC}(4)$ has 3 orbits under Kreweras complementation, one of length 8,
another of length 4, and the last of length 2.

\section{Acknowledgments}

The author thanks
Dennis Stanton, Vic Reiner, and Dennis White for explaining 
the connection with CSP.
This work was supported by the Burroughs Wellcome Fund (2005 CASI to CH)
and the National Science Foundation (DMS1815044 to CH).

\end{document}